\documentclass[12pt]{amsart}

\def\factor#1.#2.{\left. \raise 2pt\hbox{$#1$} \right/\hskip -2pt\raise -2pt\hbox{$#2$}}

\theoremstyle{plain}
\newtheorem{thm}{Theorem}[section]
\newtheorem{theorem}[thm]{Theorem}

\newtheorem{lemma}[thm]{Lemma}
\newtheorem{corollary}[thm]{Corollary}
\newtheorem{proposition}[thm]{Proposition}
\theoremstyle{definition}

\newtheorem{definition}[thm]{Definition}
\newtheorem{claim}[thm]{Claim}

\newtheorem{question}[thm]{Question}

\numberwithin{equation}{section}

\newcommand{\sC}{{\mathcal C}}

\newcommand{\sF}{{\mathcal F}}
\newcommand{\sG}{{\mathcal G}}
\newcommand{\sH}{{\mathcal H}}

\newcommand{\sN}{{\mathcal N}}
\newcommand{\sO}{{\mathcal O}}

\newcommand{\sQ}{{\mathcal Q}}

\newcommand{\sV}{{\mathcal V}}

\newcommand{\C}{{\mathbb C}}

\newcommand{\N}{{\mathbb N}}
\newcommand{\BP}{{\mathbb P}}

\newcommand{\rk}{{\rm rk}}

\title[Characteristic foliation on a hypersurface]{Characteristic foliation on a hypersurface of general type
 in  a projective symplectic manifold }
\author[Jun-Muk Hwang, Eckart Viehweg]{Jun-Muk Hwang${}^1$, Eckart Viehweg${}^2$}
\address{Jun-Muk Hwang, Korea Institute for Advanced Study, Hoegiro 87, Seoul, 130-722, Korea} \email{jmhwang@kias.re.kr}
\address{Eckart Viehweg, Universit\"at Duisburg-Essen, Mathematik, 45117 Essen, Germany} \email{viehweg@uni-due.de}
\thanks{${}^1$ supported by ASARC}
\thanks{${}^2$ supported by the DFG-Leibniz program and by the SFB/TR 45
``Periods, moduli spaces and arithmetic of algebraic varieties''}
\begin{document}

\maketitle
\section{Introduction}

In classical mechanics, it is an important question whether the
orbit of the motion of a celestial body is periodic. In the
Hamiltonian formalism, this question is formulated in terms of
symplectic geometry as follows.
 Let $(M, \omega)$ be a symplectic manifold. Given a non-singular hypersurface
$X \subset M$, the restriction of $\omega$ on the tangent space of
$X$ at  each point $x \in X$ has 1-dimensional kernel,  defining a
foliation of rank 1, which we will call the {\it characteristic
foliation of $X$} induced by $\omega$. The question on the
periodicity of orbits correspond to the following geometric
question.

\begin{question}\label{compact} Given a symplectic manifold $(M,
\omega)$ and a hypersurface $X \subset M$, when are the leaves of
the characteristic foliation of $X$ compact? \end{question}

\medskip
In Hamiltonian mechanics, $M$ corresponds to the phase space of
the mechanical system, a {\it real} symplectic manifold,   and $X$
corresponds to the level set of the energy, a {\it real}
hypersurface in $M$. It is interesting that Question \ref{compact}
makes perfect sense in the setting of complex geometry, where $M$
is a {\it holomorphic} symplectic manifold and $X$ is a {\it
complex} hypersurface. Recall that a holomorphic symplectic
manifold is a complex manifold $M$ equipped with a closed
holomorphic 2-form $\omega \in H^0(M, \Omega_M^2)$ such that
$\omega^n \in H^0(M, K_M)$ is nowhere vanishing. The holomorphic
version of Question \ref{compact} was studied in \cite{HO} and an
example where it has an affirmative answer was examined in detail.

In the introduction of \cite{HO}, the authors expected that the
answer to Question \ref{compact} would be negative for a `general'
hypersurface $X$. The aim of this paper is a verification of this
expectation when $M$ is a non-singular projective variety over
$\C$. For a holomorphic foliation of rank 1 on a projective
variety, the compactness of the leaf is equivalent to the
algebraicity of the leaf. Let us say that a holomorphic foliation
on an algebraic variety is an {\it algebraic foliation} if all of
its leaves are algebraic subvarieties. Our main result is the
following.

\begin{theorem}\label{maintheorem} Let $M$ be a non-singular projective variety of dimension $\geq 4$ with
a symplectic form $\omega$. Let $X \subset M$ be a non-singular
hypersurface of general type. Then the characteristic foliation on
$X$ induced by $\omega$ cannot be algebraic.
\end{theorem}

This applies to non-singular ample (or nef and big) hypersurfaces
$X \subset M$ because $K_M \cong \sO_M$ via the symplectic form.
Theorem \ref{maintheorem} is a direct consequence of a  more
general result:

\begin{theorem} \label{triviality} Let $X$ be a non-singular projective variety of dimension $\geq 2$
and let
$$0\longrightarrow \sQ \longrightarrow \Omega_X^1 \longrightarrow \sF \longrightarrow 0$$ be a foliation on $X$ all of whose
leaves are algebraic curves. If $\sF$ is big, i.e., its Iitaka
dimension $\kappa(\sF) = \dim(X)$,  then $\kappa(\det(\sQ)) \geq
\dim(X) -1 $.\end{theorem}

In the setting of Theorem \ref{maintheorem}, assume that the
foliation is algebraic and apply Theorem \ref{triviality}.  By the
definition of the characteristic foliation,  $\omega$ induces a
symplectic form on $\sQ$, hence $\det(\sQ) \cong \sO_X$. Since $X$
is  of general type, this implies that $\sF$ is big. This yields a
contradiction $$ 0 = \kappa ( \sO_X) = \kappa( \det(\sQ)) \geq
\dim(X) -1.$$  Thus Theorem \ref{maintheorem} is a consequence of
Theorem \ref{triviality}. Conversely, $X\subset M$ as in Theorem
\ref{maintheorem} shows that the condition that all leaves are
algebraic in Theorem \ref{triviality} is necessary.

To prove Theorem \ref{triviality}, we first develop some general
structure theory for algebraic foliations. In particular, we will
prove an \'etale version of the classical Reeb stability theorem
in foliation theory. Furthermore, for algebraic foliations by
curves, we will prove a global version of this result. Using this
general structure theorem,  the relation between Iitaka dimensions
will be obtained by borrowing a result from the theory of the
positivity of the direct image sheaves associated to families of
curves (\cite{Vi}).  The latter theory originated from the study of
the Shafarevich conjectures over function fields, and properties of
sheaves on fine moduli spaces of curves. It is amusing
to observe that these questions of modern algebraic geometry are
related to the question of the periodicity of motions of celestial
bodies.

By the decomposition theorem of \cite{Be}, non-singular projective
varieties with symplectic forms are, up to finite \'etale cover,
products of abelian varieties and projective hyperk\"ahler
manifolds. As far as we know, Theorem \ref{maintheorem} is new
even for abelian varieties. In the simplest case, we can formulate
the result explicitly as follows.

\begin{corollary} Let $A = {\C}^{2n} / \Lambda$ be an
even-dimensional principally polarized abelian variety with smooth
theta divisor. Fix any linear coordinate $(p_1, \ldots, p_n, q_1,
\ldots, q_n)$ on ${\C}^{2n}$ and let $\theta(p_1, \ldots, p_n, q_1,
\ldots, q_n)$ be the Riemann theta function on $\C^{2n}$ associated
to the period $\Lambda$. For a very general (i.e. outside a
countable union of proper subvarieties) point $(a_1, \ldots, a_n,
b_1, \ldots, b_n)$ on the theta divisor, the solution $(p_i(t),
q_i(t))$ of the Hamiltonian flow on $\C^{2n}$
$$ \frac{ d p_i}{d t} = - \frac{
\partial \theta}{\partial q_i}, \;\; \frac{ d q_i}{d t} = \frac{
\partial \theta}{\partial p_i}, \;\; i =1, \ldots, n$$ with
initial value $p_i(0) = a_i, q_i(0) = b_i, i = 1, \ldots, n,$ cannot
descend to an algebraic  curve on $A$.
\end{corollary}

It is natural to ask whether at least {\it some} leaf of the
characteristic foliation in Theorem \ref{maintheorem} can be an
algebraic curve.  This question is completely out of the reach of
the method employed in the current paper. We cannot even make a
guess whether the answer would be affirmative or not.

\section{\'Etale Reeb stability for algebraic foliations}

Let $X$ be a non-singular projective variety over $\C$. A foliation
of $X$ is given by an exact sequence of locally free sheaves $$ 0
\longrightarrow \sQ \longrightarrow \Omega^1_X \longrightarrow \sF
\longrightarrow 0$$ where $\sQ$ is integrable, i.e., through each
point $x \in X$, there exists a complex submanifold $C$ such that
$\sQ$ corresponds to the conormal bundle of $C$ at every point of
$C$. This submanifold $C$ is called the leaf of the foliation
through $x$.  We say that the foliation is  algebraic if each leaf
is an algebraic subvariety of $X$. Our aim in this section is to
describe the behavior of the leaves of an algebraic foliation as a
family of algebraic subvarieties. For that purpose, we need to
recall some standard results on the structure of differentiable
foliations.

Let $X$ be a differentiable manifold with a differentiable
foliation. A {\it transversal section} at a point $x \in X$ means
a (not necessarily closed) submanifold $S$ through $x$ whose
dimension is equal to the codimension of the leaves such that the
intersection of each leaf of the foliation with $S$ is
transversal (or empty). Let $C$ be the leaf through $x$. A choice
of a transversal section $S$ at $x$ determines a group
homomorphism
$$ \pi_1(C, x) \longrightarrow {\rm Diff}_x(S),$$ called the {\it holonomy
homomorphism}, from the fundamental group of $C$ to the group ${\rm
Diff}_x(S)$ of germs of the diffeomorphisms of $S$ at $x$. For a
precise definition of this homomorphism, we refer the readers to
\cite[Section 2.1]{MM}. Roughly speaking, a loop $\gamma$ on $C$
representing an element of $\pi_1(C,x)$ acts on $S$ by moving a
point $y \in S$ close to $x$ along the leaf through $y$ following
$\gamma$. The image of the holonomy homomorphism will be called the
{\it holonomy group} of the leaf $C$. The isomorphism class of this
group depends only on $C$, independent from the choice of $x$ and
$S$. The following is a well-known criterion for the finiteness of
the holonomy group.

\begin{proposition}\label{finitehol} Let $X$ be a differentiable manifold with a
foliation all of whose leaves are compact. The holonomy group of a
leaf $C$ is finite if  there exist  a transversal section $S$ at a
point $x \in C$ and a fixed positive integer $N$ such that the
cardinality of the intersection of $S$ with any leaf of the
foliation is bounded by $N$. \end{proposition}

\begin{proof} This is contained in \cite[Theorem 4.2]{Ep} and the proof can be
found in \cite[Section 7]{Ep}. In fact, all we really need is the
simple fact that if a finitely generated group $G$ acts
effectively on a set $S$ such that the cardinality of each orbit
is bounded by a positive integer $N$, then the group is finite. We
recall the proof for the reader's convenience. Denote by ${\rm
S}_r$ the permutation group of $r$ points. Since $G$ is finitely
generated there are only finitely many homomorphisms $G \to {\rm
S}_r$. Let $H \subset G$ be the intersection of the kernels of all
such homomorphisms for $r \leq N$. Then $G/H$ is finite. Since
each orbit of $G$ in $S$ determines a group homomorphism $G \to
{\rm S}_r$ for some $r \leq N$, $H$ must act trivially on $S$.
Thus $H$ is the trivial subgroup. It follows that $G$ is finite.
\end{proof}

We recall the construction of the {\it  flat bundle foliation} in
\cite[p.17]{MM}. Let $G$ be a finite group which acts freely  on a
manifold $\tilde{C}$ on the right.  Suppose $G$ acts effectively on
another manifold $S$  on the left with a fixed point $x \in S$. Let
$\tilde{C} \times_G S$ be the quotient of $\tilde{C} \times S$ by
the equivalence relation $(y g, s) \sim (y, gs)$ for $g \in G$ and
$(y, s) \in \tilde{C} \times S$. Let $ C \subset \tilde{C} \times_G
S$ be the image of $\tilde{C} \times \{x\}$.
 We have the commutative diagram
\begin{equation}\label{com} \begin{array}{ccc} \tilde{C} \times S & \longrightarrow &
\tilde{C} \times_G S \\ \downarrow & & \downarrow \\ S &
\longrightarrow & G\backslash S. \end{array} \end{equation} The
foliation on the manifold $\tilde{C} \times_G S$ given by the
vertical fibers is called the {\it flat bundle foliation} arising
from the actions of $G$ on $\tilde{C}$ and $S$.  $C$ is a leaf of
this foliation and $G$ is the holonomy group of the leaf $C$. For
any $y \in \tilde{C}$, the image of $\{y \} \times S$ gives a
transversal section of this foliation at the image $\bar{y} \in C$.
Then  $\tilde{C}$ is the $G$-Galois cover of $C$ associated to the
holonomy homomorphism $\pi_1(C, \bar{y}) \to G$. The following is
easy to check.

\begin{lemma}\label{normalization} Let $S' \subset \tilde{C} \times_G S$ be a closed
submanifold  with $S' \cap C =: \{y\}$ which is a
 transversal section of the flat bundle foliation at $y \in C$.  Then \eqref{com} factors through
 the (set-theoretic) fiber product of $S' \to G\backslash S$ and $\tilde{C} \times_G S \to G \backslash S$
via a finite map $$ \tilde{C} \times S  \longrightarrow S'
\times_{G \backslash S} (\tilde{C} \times_G S)$$ which is
one-to-one over a dense open subset. \end{lemma}

We say that a subset of a manifold with a foliation is {\it
saturated} if it is the union of leaves intersecting it. The
following is the classical Reeb stability theorem whose proof can
be found in \cite[Theorem 2.9]{MM}.

\begin{theorem}\label{Reeb}   For a  differentiable
manifold $X$ with a foliation, suppose  $C$ is a compact leaf with
finite holonomy group $G$. Then there exist a saturated open
neighborhood $U$ of $C$ in $X$ and a transversal section $S$ in
$U$ such that denoting by $G$ the finite holonomy group acting on
$S$ and by $\tilde{C} \to C$ the $G$-Galois covering, there exists
a diffeomorphism $\tilde{C} \times_G S \cong U$ such that the
foliation on $U$ correspond to the flat bundle foliation on
$\tilde{C} \times_G S$. \end{theorem}

Now assume that $X$ is a complex manifold with  a holomorphic
foliation. Then we have the following holomorphic version of Reeb
stability.

\begin{theorem}\label{holReeb}  For  a complex manifold $X$ with a
holomorphic foliation, suppose that $C$ is a compact leaf with
finite holonomy group $G$.  Then there exist a saturated open
neighborhood $U$ of $C$ in $X$,  a holomorphic transversal section
$S$ in $U$ with a $G$-action,  an unramified $G$-Galois cover
$\tilde{U} \to U$,  a smooth proper morphism $h: \tilde{U} \to S$
and a proper morphism $g: U \to G\backslash S$ satisfying the
commutative diagram \begin{equation}\label{holcom}
\begin{array}{ccc}
\tilde{U}  & \longrightarrow & U \\
h \downarrow & & \downarrow g \\ S & \longrightarrow & G\backslash
S.
\end{array} \end{equation}
 Moreover, for each closed
submanifold $\Sigma \subset U$ which intersects all leaves
transversally, the normalization of the fiber product $\Sigma
\times_{G\backslash S} U$ is an unramified cover of $U$.
\end{theorem}

\begin{proof} We can apply Theorem \ref{Reeb}.  The diagram \eqref{holcom}
is just \eqref{com} where $U= \tilde{C} \times_G S, \tilde{U} =
\tilde{C} \times S$ and  $\tilde{U}$ is given the complex
structure as an unramified  covering of $U$. The differentiable
maps $h$ and $g$ are holomorphic maps because the foliation is
holomorphic. The last statement is a consequence of Lemma
\ref{normalization}.
\end{proof}

We can apply this theorem to algebraic foliations because:

\begin{proposition}\label{chow} Let $X$ be a non-singular projective variety
with an algebraic foliation. Then each leaf has finite holonomy
group. Denoting by ${\rm Chow}_X$ the Chow variety of $X$, there
exists a natural morphism $\mu: X \to {\rm Chow}_X$ sending all
points on a leaf $C$ with holonomy group $G_C$ to the cycle $|G_C|
\cdot C$.\end{proposition}

\begin{proof}
Given a leaf $C$ and a point $x \in C$, we can find a complete
intersection $S$ of very ample hypersurfaces which intersects $C$
transversally with $x \in C \cap S$. In an analytic neighborhood
of $C$, a component of $S \cap U $ is a transversal section. By
Noetherian induction, the intersection number of each leaf of the
foliation with $S$ is bounded by a positive number $N$. Thus we
can apply Proposition \ref{finitehol} to conclude that the
holonomy group is finite. The fact that the cycles $|G_C| \cdot C$
form a nice family follows from the local description of the
family of leaves in Theorem \ref{holReeb}.
\end{proof}

For the next theorem, we  need the following lemma.

\begin{lemma}\label{lem2} In the setting of Proposition \ref{chow},
  let $\Sigma \subset X$ be a subvariety such
that $$\nu:= \mu|_{\Sigma}: \Sigma \longrightarrow \mu(X)$$ is  a finite
morphism. Let $M_0 \subset \mu(X)$ be a connected analytic open
subset such that for each point $y \in M_0$, the reduction of the
fiber $\mu^{-1}(y)_{\rm red}$ intersects $\Sigma$ transversally. Set
$X_0 := \mu^{-1}(M_0)$ and $\Sigma_0 := \nu^{-1}(M_0)$.  Then the
normalization of the fiber product of $\mu_0: X_0 \rightarrow M_0$
and $\nu_0: \Sigma_0 \rightarrow M_0$ is an unramified covering of
$X_0$.
\end{lemma}

\begin{proof} The statement is local on $\mu(X)$. So we can verify it
for any neighborhood of a given point $y \in \mu(X)$. In other
words, we may assume that $X_0$ is contained in the neighborhood $U$
of Theorem \ref{holReeb}. Then this is immediate from the last
statement in Theorem \ref{holReeb}.  \end{proof}

The following is the \'etale version of local Reeb stability
theorem for an algebraic foliation.

\begin{theorem}\label{etaleReeb}  Let $X$ be a non-singular projective variety with
an algebraic foliation \begin{equation}\label{eq.1} 0\longrightarrow \sQ \longrightarrow
\Omega_X^1 \longrightarrow \sF \longrightarrow 0.
\end{equation}  Then
for each leaf $ C \subset X$, there exists an \'etale neighborhood
$\tau: U \to X$ of $C$, a smooth projective morphism $h : U \to M$
and isomorphisms $$\tau^* \sQ \cong h^* \Omega^1_M \mbox{ and }
\tau^* \sF \cong \Omega^1_{U/M}$$ such that the pullback of
\eqref{eq.1} is isomorphic to the tautological exact sequence $$
0\longrightarrow h^*\Omega_M^1 \longrightarrow \Omega_U^1 \longrightarrow \Omega^1_{U/M} \longrightarrow 0. $$
\end{theorem}

\begin{proof} Just take a general complete intersection $\Sigma \subset X$
of very ample divisors  intersecting $C$ transversally. Then apply
Lemma \ref{lem2}, with $M_0$ the Zariski open subset where the
reduction of fibers of $\mu$ intersect $\Sigma$ transversally. The
\'etale neighborhood $U$ is given by the normalization of the
fiber product $X_0 \times_{M_0} \Sigma_0$. The existence of the
smooth morphism $U \to M$ follows from Theorem \ref{holReeb}.
\end{proof}

\section{Global \'etale Reeb stability for algebraic foliations by curves}

For algebraic foliations by curves, we can globalize Theorem
\ref{etaleReeb}. The essential point is the existence of the moduli scheme
$M_g$ of curves of genus $g$, or of the moduli schemes $M_g^{[N]}$
of curves of genus $g$ with a level $N$ structure, which are fine for $N\geq 3$.

\begin{lemma}\label{fine}
Let $f:V\to W$ be a smooth morphism of curves and $N\in \N$. Then
there exists an \'etale finite morphism $\tilde{W}\to W$ such that
$$
\tilde{V}=V\times_W\tilde{W} \longrightarrow \tilde{W}
$$
carries a level
$N$-structure.  In particular, for $N\geq 3$ the family
$\tilde{V}\to \tilde{W}$ is the pullback of the universal family
over the fine moduli scheme $M_g^{[N]}$ of curves with a level $N$
structure.
\end{lemma}
\begin{proof}
One can choose a level $N$ structure if  the $N$-division points
of the relative Jacobian are generated by sections. Since the
$N$-division points of a family of abelian varieties are \'etale
and finite over the base, this can be achieved over a finite
\'etale cover.
\end{proof}

The global version of Theorem \ref{etaleReeb} is the following.

\begin{theorem}\label{morphism}
Let $X$ be a non-singular projective variety with an algebraic
foliation \begin{equation}\label{eq.m} 0\longrightarrow \sQ \longrightarrow \Omega_X^1 \longrightarrow
\sF \longrightarrow 0.
\end{equation} whose leaves are curves of genus $\geq 2$. Then there
exist a generically finite projective morphism $\sigma:V \to X$, a
non-singular projective variety $W$, a smooth projective morphism
$f: V\to W$ and injections
$$
\sigma^*\sQ \stackrel{\alpha}{\longrightarrow} f^*\Omega^1_W, \ \
\mbox{ and } \ \  \sigma^*\sF \stackrel{\cong}{\longrightarrow}
\Omega^1_{V/W},
$$
such that the pullback of \eqref{eq.m} is a subcomplex of the
tautological exact sequence
$$
0\longrightarrow f^*\Omega_W^1 \longrightarrow \Omega_V^1
\longrightarrow \Omega^1_{V/W} \longrightarrow 0. $$ For each point
$w\in W$ one can choose a neighborhood $W_0$ and an \'etale open
neighborhood $U$ of the image of $f^{-1}(w)$ in $X$, satisfying the
condition in Theorem \ref{etaleReeb}. In particular $f:V\to W$ is a
smooth family of curves. Moreover, we can assume that the associated
classifying morphism $W\to M_g$ factors like
$$
W \stackrel{\varphi'}{\longrightarrow} M_g^{[N]}
$$ for a given positive integer $N$.
 \end{theorem}

\begin{proof}
Let us choose a finite set of \'etale neighborhoods appearing in
Theorem \ref{etaleReeb}, say $\tau_i:U_i\to X$ for $i\in
\{1,\ldots,\ell\}$, such that
\begin{equation}\label{eq.3}
\bigcup_{i=1}^\ell \tau_i(U_i)= X.
\end{equation}
 The families $h_i:U_i\to M_i$ induce morphisms $\phi_i:M_i\to
{\rm Chow}_X$ such that
$$ \begin{array}{ccc} U_i & \stackrel{\tau_i}{\longrightarrow} & X
\\ h_i \downarrow & & \downarrow \mu \\ M_i &
\stackrel{\phi_i}{\longrightarrow} & {\rm Chow}_X, \end{array} $$
where $\mu$ is the morphism in Proposition \ref{chow}.

Next we fix some projective compactification $\bar{M}_i$ such that
$\phi_i$  extends  to a morphism $\bar{\phi}_i:\bar{M}_i\to {\rm
Chow}_X$. Replacing $\bar{M}_i$ by the Stein factorization we may
as well assume that $\bar{\phi}_i$ is finite. Let $\bar{M}'$ be an
irreducible component of
$$
\bar{M}_1 \times _{{\rm Chow}_X} \cdots \times_{{\rm Chow}_X}
\bar{M}_\ell,
$$
with induced morphism $\bar{\phi}':\bar{M}' \to {\rm Chow}_X$. Let
$W$ be the normalization of $\bar{M}'$ in the Galois hull of the
function field $\C(\bar{M}')$ over $\C(\bar{\phi}'(\bar{M}'))$.
Hence writing $\bar{\phi}:W \to {\rm Chow}_X$  for the induced
morphism, there is a finite group $G$ acting on $\bar{M}$ with
quotient $\bar{\phi}(W)$. The condition \eqref{eq.3} implies that
$$
\bigcup_{i=1}^\ell \phi_i(M_i) = \bar{\phi}(W) = \mu(X).$$

Let $\tilde{M}_i$ denote the preimage of $M_i$ under
$$
\bar{M}' \subset \bar{M}_1 \times _{{\rm Chow}_X} \cdots
\times_{{\rm Chow}_X} \bar{M}_\ell \stackrel{{\rm
pr}_i}{\longrightarrow} \bar{M}_i.
$$
By pullback there is a smooth projective morphism $\tilde{U}_i \to
\tilde{M}_i$. For all $\gamma\in G$ one obtains the pullback
family
$$
\tilde{U}^\gamma_i \longrightarrow \tilde{M}^\gamma_i :=
\gamma^{-1}(\tilde{M}_i).
$$
For different $i, \ i'$ and for $\gamma$ and $\gamma'$ the closed
fibres of those families coincide on
$$
\tilde{M}^\gamma_i \cup \tilde{M}^{\gamma'}_{i'}.
$$
In fact, the isomorphism class of a fibre is determined by the
image in $X$, hence it is invariant under $G$ and independent of
the \'etale neighborhood.

In particular, the morphisms $\tilde{U}^\gamma_i \to M_g$
mapping a point $w$ to the moduli point of the isomorphism class
of the fibre over $w$ glue to a morphism $W\to M_g$. Replacing
$W$ by a finite covering, we may assume by Lemma \ref{fine} that
this morphism factors through the fine moduli scheme
$M_g^{[N]}$. Then the different families over the open subsets
$\tilde{M}^\gamma_i$ are pullbacks of the universal family over
$M_g^{[N]}$. hence they coincide over the two by two
intersections, and glue to a smooth family $V\to W$.

By abuse of notations, we replace $W$ by a desingularization and
$f:V\to W$ by the pullback family. It satisfies all the required
properties. \end{proof}

\section{Positivity property of algebraic foliations by curves}

Let us recall some notions of  positivity for locally free sheaves.
\begin{definition}\label{bigdef}
Let $\sG$ be a locally free sheaf on a quasi-projective variety $Z$
and let $Z_0 \subset Z$ be an open dense subvariety. Let $\sH$ be an
ample invertible sheaf on $Z$.
\begin{enumerate} \item[1.]
$\sG$ is {\em globally generated over $Z_0$} if the natural morphism
$$
H^0 (Z, \sG) \otimes \sO_Z \longrightarrow \sG
$$
is surjective over $Z_0$.
\item[2.] $\sG$ is {\em ample with respect to $Z_0$} if for some $k
>0$, $S^{k}(\sG) \otimes \sH^{-1}$ is globally generated over
$Z_0$. In particular, $\sG$ is {\em ample} if it is ample with
respect to $Z_0=Z$.
\item[3.] $\sG$ is {\em big} if it is ample with respect to some
open dense subvariety $Z_0$. If $\rk(\sG) =1$, this is equivalent to
saying that the Iitaka dimension $\kappa(\sG) = \dim Z$.
\end{enumerate}
\end{definition}
In the literature one finds a second definition for bigness of a
locally free sheaf, requiring $\sO_{\BP(\sG)}(1)$ to be big on the
projective bundle $\pi:\BP(\sG) \to Z$ induced by $\sG$. Our notion
of bigness is stronger. It is equivalent to the ampleness of
$\sO_{\BP(\sG)}(1)$ with respect to an open set of the form
$\pi^{-1}(Z_0)$.

\begin{lemma}\label{biglemma}  Let $\sG$ be a  locally free sheaf on a quasi-projective
non-singular variety $Z$.
\begin{enumerate} \item[(i)] If $\sG$ is big, then for a locally free sheaf $\sG'$ and
a non-zero  homomorphism $\eta: \sG \to \sG'$, the sheaf
$\det(\eta(\sG))$ is big. \item[(ii)] If $\sG$ is ample, then for
any generically finite morphism $\rho: Y \to Z$, the pullback
$\rho^*\sG$ is big.
\end{enumerate}
\end{lemma}
The proof of the lemma is straight-forward. Let us just point out, that we define $\det(\eta(\sG))$ to be  $\iota_*\det(\eta(\sG|_{Z'}))$, where $\iota:Z'\to Z$ is the largest open subscheme with $\eta(\sG|_{Z'})$ locally free.
The assumption ``non-singular'' is needed to get an invertible sheaf. Without it one would have to allow in Definition \ref{bigdef} torsion free coherent sheaves $\sG$, making the notations more complicated.

We will need the following result  from \cite[Proposition 2.4]{Vi}.

\begin{proposition}\label{moduli} Let $ \mathfrak{f}:
\sC_g \to M_g^{[N]}, N \geq 3,$ be the universal family over the
fine moduli scheme of curves with level $N$-structures.  Then
$\mathfrak{f}_*\omega_{\sC_g/M_g^{[N]}}^\nu$ is  ample  for all
$\nu\geq 2$.
\end{proposition}

Theorem \ref{triviality} is a direct consequence of the following.

\begin{proposition}\label{subsheaf}  In the setting of Theorem \ref{morphism}, let $v =
\kappa(\sF) -1$. Then one has:
\begin{enumerate}
\item[a)] $v ={\rm Var}(f)$,  the dimension of the image of $W$ in
$M_g$. \item[b)] There is a subsheaf $\sV\subset \sigma^*\sQ$ of
rank $v$ with $\kappa(\det(\sV))=v$.
\end{enumerate}
\end{proposition}

\begin{proof}
Fix $N \geq 3$. Since $M_g^{[N]}$ is a fine moduli scheme, $V\to
W$ is the pullback of the universal family $\sC_g\to M_g^{[N]}$
under the morphism $\varphi':W\to M_g^{[N]}$. Consider a
factorization
$$
W \stackrel{\varphi}{\longrightarrow} Z \stackrel{\rho}{\longrightarrow} M_g^{[N]}.
$$
with $\varphi$ surjective and with connected fibres, and with $\rho$ generically finite.
Blowing up, and replacing the families with the pullbacks, we may assume that $Z$ is non-singular.
Let us write $g:T\to Z$ for the pullback of the universal family to $Z$. Then
$V\cong T\times_Z W$ and the second projection defines a morphism $p:V\to T$.

Since $\omega_{V/W}=p^*\omega_{T/Z}$ one finds that
$$
\kappa(\sF)=\kappa(\sigma^*
\sF)=\kappa(\omega_{V/W})=\kappa(\omega_{T/Z}).
$$

\begin{claim}\label{claim1}
The invertible sheaf $\omega_{T/Z}$ is big.
\end{claim}
Since $\dim(Z)= {\rm Var}(f)$, Claim \ref{claim1} shows $$v+1=
\kappa(\omega_{T/Z})=\dim(Z)+1= {\rm Var}(f), $$ proving a).

\begin{proof}[Proof of Claim \ref{claim1}]
From Proposition \ref{moduli},
$\mathfrak{f}_*\omega^2_{\sC_g/M_g^{[N]}}$ is ample. Let $\sN$ be
an ample invertible sheaf on $M_g^{[N]}$. For some $k \gg 1$ the
sheaf $\sN^{-1} \otimes
S^k(\mathfrak{f}_*\omega^2_{\sC_g/M_g^{[N]}})$ is globally
generated. Writing $\sH$ for the pullback of $\sN$ to $Z$ one finds
by base change that $\sH^{-1} \otimes S^k(g_*\omega^2_{T/Z})$ is
globally generated. Using the multiplication map one gets an
inclusion
$$
\bigoplus \sH^\ell \hookrightarrow g_* \omega^{2\cdot \ell \cdot
k}_{T/Z}
$$ where $\bigoplus \sH^\ell$ denotes the sum of $\rk(g_* \omega^{2\cdot \ell \cdot k}_{T/Z})$ copies of $\sH^\ell$. So $h^0(T,\omega^{2\cdot \ell \cdot
k}_{T/Z})$ is larger than a polynomial in $\ell$ of degree
$\dim(Z)+1$ and with positive leading coefficient.
\end{proof}
To prove b), recall that we have the Kodaira-Spencer homomorphism
$$\mathfrak{f}_* \omega^2_{\sC_g/M_g^{[N]}} \longrightarrow
\Omega^1_{M_g^{[N]}}$$ whose pullback $$\eta: g_*\omega^2_{T/Z}
\longrightarrow \Omega^1_Z$$ must be surjective over a Zariski open subset of
$Z$ because $Z$ is generically finite over $M_g^{[N]}$. Let
$\Omega \subset \Omega^1_Z$ be the image of the homomorphism $\eta$.
Since $\mathfrak{f}_* \omega^2_{\sC_g/M_g^{[N]}}$ is ample,
$g_*\omega^2_{T/Z} = \rho^* \mathfrak{f}_*
\omega^2_{\sC_g/M_g^{[N]}}  $ is big by Lemma \ref{biglemma} (ii),
hence $\det(\Omega)$ is big by Lemma \ref{biglemma} (i). Thus
$$ \kappa( \varphi^*\det(\Omega)) = \kappa(\det( \varphi^*
\Omega)) = \kappa(\det(\Omega)) = v.$$

Let $\sV:= f^*\varphi^*\Omega$. Since
$\det(\sV)=f^*\varphi^*\det(\Omega)$, its Iitaka dimension is
equal to $\dim(Z)$. Thus to prove b), it suffices to verify:
\begin{claim}\label{subsheaf2}
$\sV = f^*\varphi^*\Omega$ is a subsheaf of $\sigma^*\sQ$.
\end{claim}
{\em Proof of Claim \ref{subsheaf2}.} One has a natural inclusions
$\varphi^*\Omega^1_Z \to \Omega_W^1$, and hence $f^*\varphi^*\Omega
\to f^*\Omega_W^1$. To see that its image lies in the smaller sheaf
$\sigma^*\sQ$ is a local question. So it will be sufficient to
verify this in the neighborhood $W_0$ considered in Theorem
\ref{morphism}. Over $W_0$ the morphism is the pullback of the
morphism $h: U\to M$ in Theorem \ref{etaleReeb}. So the morphism
$W_0\to M_g^{[N]}$ factors like $W_0\to M \to M_g^{[N]}$ and the
pullback of $\mathfrak{f}_* \omega^2_{\sC_g/M^{[N]}_g}$ is $h_*
\omega^2_{U/M}$, which is sent to $\Omega^1_M$ by the
Kodaira-Spencer map. It follows that the pullback of $\Omega$ lies
in  $h^*\Omega^1_{M}=\tau^*\sQ$.
\end{proof}

\begin{proof}[Proof of Theorem \ref{triviality}]   Since
$\sF$ is big,  general leaves have genus $\geq 2$. By Theorem
\ref{Reeb}, every leaf has genus $\geq 2$. Moreover, the number
$v$ in Proposition \ref{subsheaf}  is equal to $\dim(X)-1$, hence
to ${\rm rank}(\sQ)$. So the subsheaf $\sV$ in Proposition
\ref{subsheaf} has the same rank as $\sigma^*\sQ$ and
$$\kappa(\det(\sQ)) \geq \kappa(\det(\sV)) = v.$$ \end{proof}

\end{document}